\documentclass[10pt]{amsart}
\usepackage{amsmath,amssymb,latexsym}

\newtheorem{conj}[equation]{Conjecture}
\newtheorem{prop}[equation]{Proposition}
\newtheorem{dfn}[equation]{Definition}

\newtheorem{thm}[equation]{Theorem}
\newtheorem{cor}[equation]{Corollary}
\newtheorem{lem}[equation]{Lemma}

\theoremstyle{remark}
\newtheorem{rem}[equation]{Remark}

\makeatletter\@addtoreset{equation}{section}\makeatother

\def\B{{\mathcal{B}}}
\def\E{{\mathcal{E}}}

\def\C{{\mathcal{C}}}

\def\Q{{\mathbb{Q}}}
\def\Z{{\mathbb{Z}}}

\def\e{{\varepsilon}}
\def\h{{\mathfrak{h}}}

\def\YD{{\mathcal{YD}}}
\def\af{{\rm{aff}}}

\newcommand{\arr}{\begin{picture}(35,5) \put(3,2.5){\vector(1,0){29}} \end{picture}}
\newcommand{\arrr}{\begin{picture}(53,5) \put(3,2.5){\vector(1,0){47}} \end{picture}}
\newcommand{\arrl}[1]{\stackrel{#1}{\arr}}
\newcommand{\arrll}[1]{\stackrel{#1}{\arrr}}

\title[Nichols-Woronowicz model of the $K$-theory of generalized flag varieties]
{Alcove path and Nichols-Woronowicz model of the equivariant $K$-theory of generalized flag varieties}

\author{Cristian Lenart}
\address{Department of Mathematics and Statistics, 
State University of New York at Albany, Albany, NY 12222, USA}
\email{lenart@albany.edu}
\urladdr{http://math.albany.edu/math/pers/lenart}
\author{Toshiaki Maeno}
\address{Department of Electrical Engineering, Kyoto University, Sakyo-ku, Kyoto 606-8501, Japan}
\email{maeno@kuee.kyoto-u.ac.jp}

\thanks{Cristian Lenart was supported by National Science Foundation grant DMS-0403029.\\
\indent Toshiaki Maeno was supported by Grant-in-Aid for Scientific Research.}

\begin{document}
\begin{abstract}
Fomin and Kirillov initiated a line of research into the realization of the 
cohomology and $K$-theory of generalized flag varieties $G/B$ as commutative 
subalgebras of certain noncommutative algebras.  This approach has 
several advantages, which we discuss. This paper contains the most comprehensive 
result in a series of papers related to the mentioned line of research. 
More precisely, we give a model for the $T$-equivariant $K$-theory 
of a generalized flag variety $K_T(G/B)$ in terms of a certain braided 
Hopf algebra called the Nichols-Woronowicz 
algebra. Our model is based on the Chevalley-type multiplication formula 
for $K_T(G/B)$ due to the first author and Postnikov; this formula is stated 
using certain operators defined in terms of so-called alcove paths (and 
the corresponding affine Weyl group). Our model is derived using a 
type-independent and concise approach. 
\end{abstract}

\maketitle

\begin{center} {\small\it{
Dedicated to Professor Kenji Ueno on the occasion of his sixtieth birthday}}
\end{center}
\section*{Introduction}

Schubert calculus evolved from the calculus of enumerative geometry to the 
study of geometric, algebraic, and combinatorial aspects related to various 
cohomology settings for algebraic homogeneous spaces. An important problem 
in Schubert calculus is to find a combinatorial description of the 
cohomology of generalized flag varieties $G/B$ 
(where $G$ is a semisimple Lie group and 
$B$ a Borel subgroup); more general algebras were also considered, such 
as $K$-theory and quantum cohomology. 

Fomin and Kirillov \cite{FK} initiated a line of research into the realization 
of the mentioned algebras as commutative subalgebras of certain noncommutative 
algebras. The idea is to map certain cohomology classes to elements of 
the noncommutative algebra which are thought of as multiplication operators 
acting on cohomology. This approach has several advantages. 
First of all, it allows us to recursively construct the basis of Schubert 
classes by using a bottom-up approach, rather than the traditional top-down 
approach based on {\em divided difference operators} (see \cite{FK}). 
Secondly, given that these constructions have certain multiplication operators 
(related to the {\em Chevalley multiplication} formula) built into 
them, they are readily amenable to deriving more general multiplication 
formulas. This approach was succesfully used in \cite{Po}, where a {\em 
Pieri-type multiplication} formula in the cohomology and quantum cohomology 
of the flag variety $Fl_n$ of type $A$ was derived using the Fomin-Kirillov 
construction. Thirdly, in certain cases, the divided difference operators 
in cohomology correspond to some natural operators acting on the noncommutative 
algebra (e.g., certain twisted derivations in the case of the 
{Nichols-Woronowicz algebra}, see below).

We will now give more details about the line of research initiated by Fomin 
and Kirillov. They constructed a combinatorial 
model for the cohomology of the flag variety ${F}l_n$ 
as a subalgebra of a certain noncommutative algebra $\E_n$ defined by 
quadratic relations. More precisely, $\E_n$ is an algebra generated over 
${\bf Z}$ by the generators $[i,j],$ $1\leq i\ne j \leq n,$ subject to 
the following relations: 
\begin{enumerate}
\item[(0)] $[i,j]=-[j,i],$ 
\item[(1)] $[i,j]^2=0,$ 
\item[(2)] $[i,j][k,l]=[k,l][i,j],$ if $\{i,j \}\cap \{ k,l \} =\emptyset,$ 
\item[(3)] $[i,j][j,k]+[j,k][k,i]+[k,i][i,j]=0.$
\end{enumerate}
Fomin and Kirillov defined the commuting family of {\em Dunkl elements} 
$\theta_1,\ldots,\theta_n$ in $\E_n$ by 
\[ \theta_i:=\sum_{j\not=i}[i,j], \]
and proved that the subalgebra generated by 
them is isomorphic to the cohomology ring of ${F}l_n.$ 
It is remarkable that the algebra $\E_n$ has a natural quantum deformation. 
The deformed algebra $\widetilde{\E}_n$ is an algebra over the polynomial 
ring 
$\Z[q_1,\ldots ,q_{n-1}]$ obtained by replacing relation $(1)$ above by 
the relation 
\begin{enumerate}
\item[(1$'$)] $[i,j]^2=0,$ for $j>i+1,$ and $[i,i+1]^2=q_i.$ 
\end{enumerate}
It was conjectured in \cite{FK} that the subalgebra of $\widetilde{\E}_n$ 
generated by the Dunkl elements is isomorphic to the quantum cohomology 
ring $QH^*({F}l_n).$ This conjecture was proved by Postnikov \cite{Po}. 
A generalization to other root systems of the Fomin-Kirillov quadratic 
algebra and of the realization of $H^*(Fl_n)$ inside it is given in \cite{KM1}. 
A realization of the $K$-theory of $Fl_n$ as a commutative subalgebra 
of $\E_n$ was conjectured in \cite{Le,LY} and proved in \cite{KM3}. This 
realization is based on the so-called multiplicative Dunkl elements in 
\cite{Le,LY,Ru}; the latter were proved, in \cite{Le1}, to describe the 
Chevalley-type (or Monk-type) multiplication formula for {\em Grothendieck 
polynomials} (which are polynomial representatives for Schubert classes 
in $K(Fl_n)$).

The Hopf algebra structure related to the algebra $\E_n$ was studied 
by Fomin and Procesi \cite{FP}. The relationship between the algebra $\E_n$ 
and a braided Hopf algebra called the {\em Nichols-Woronowicz algebra} was 
pointed out by Milinski and Schneider \cite{MS}. 
Conjecturally, the algebra $\E_n$ is isomorphic to the Nichols-Woronowicz 
algebra $\B(V_{S_n})$ associated to a {\em Yetter-Drinfeld module} $V_{S_n}$ 
over the symmetric group $S_n.$ More generally, given an arbitary finite 
Coxeter group $W$, one can consider the Nichols-Woronowicz algebra $\B(V_W)$ 
associated to the {Yetter-Drinfeld module} $V_W$ over $W$. 

Bazlov \cite{Ba} constructed a model for the coinvariant algebra of a finite 
Coxeter group $W$ as a subalgebra of the corresponding Nichols-Woronowicz 
algebra $\B(V_W).$ The natural braided differential operators acting on the 
Nichols-Woronowicz algebra play a key role in his argument. 
In \cite{KM2}, Kirillov and the second author constructed a model for the 
small quantum cohomology ring $QH^*(G/B)$ in terms of the Nichols-Woronowicz 
algebra. 
In \cite{KM3}, Kirillov and the second author constructed a model for $K(G/B)$ 
in $\B(V_W)$ when the Lie group $G$ is of classical type or of type $G_2.$

The purpose of this paper is to construct a model for the $T$-equivariant 
$K$-theory of an 
arbitrary generalized flag variety $K_T(G/B)$ as a subalgebra of the 
Nichols-Woronowicz algebra $\B(V_W)$ 
associated to the Yetter-Drinfeld module $V_W$ over the corresponding Weyl 
group $W.$ Our model is based on the Chevalley-type multiplication formula 
for $K_T(G/B)$ in \cite{LP}; this formula is stated using certain operators 
defined in terms of so-called {\em alcove paths} (and the corresponding 
affine Weyl group). We define the {\em multiplicative Dunkl elements} 
in $\B(V_W)$ based on the path operators in \cite{LP}, 
and show that the commutative subalgebra generated by them is isomorphic 
to the $K$-ring $K(G/B);$ this result is then extended to equivariant $K$-theory. 
Hence, the present paper extends the results in \cite{Ba,FK,KM1,KM3}, and 
contains the most comprehensive result related to the cohomology and $K$-theory 
of $G/B$ in the line of research initiated by Fomin and Kirillov. 

Being based on the Nichols-Woronowicz algebra and the path operators in 
\cite{LP}, our model uses 
a type-independent approach, which also has the advantage of being concise. 
The quantization of our model, which is expected to give a model for the 
{\em quantum $K$-ring} of the flag variety $G/B$ \cite{GL,L}, will be 
discussed in \cite{KM4}. 

We now offer a preview of the path operators on which our construction is 
based. Let $P$ be the weight lattice in the Cartan subalgebra $\h$ of the Lie algebra 
of the simple Lie group $G.$ Denote by $T$ the maximal torus in $G$, and 
by $R(T)$ the representation ring of $T.$ 
The $K$-ring $K(G/B)$ has a presentation as a quotient of the group algebra 
of the weight lattice $P;$ a similar description can be given for the 
$T$-equivariant $K$-theory $K_T(G/B)$, as follows.
Let us consider an algebra isomorphism 
$\iota : \Z [P] = \bigoplus_{\lambda \in P} \Z \cdot e^{\lambda}
\rightarrow R(T)$ such that $\iota (e^{\lambda})$ 
is the character $\chi^{\lambda}$ corresponding to the weight $\lambda.$ 

The $T$-equivariant $K$-algebra $K_T(G/B)$ is isomorphic to the quotient algebra 
$R(T)\otimes \Z [P]/J,$ where 
$J$ is the ideal $(1\otimes f-\iota (f) \otimes 1,\; \; f\in \Z [P]^{W}).$ 
The first author and Postnikov \cite{LP} introduced the path operator $R^{[\lambda]}$ 
acting on $K_T(G/B)$ in order to derive a Chevalley-type formula in $K_T(G/B)$ which 
describes the multiplication by the class of the line bundle ${\mathcal L}_{\lambda}$ 
on $G/B$ associated to the weight $\lambda.$ 
The path operator $R^{[\lambda]}$ is defined by using the alcove path 
$A_0 \stackrel{-\beta_1}{\longrightarrow} \cdots 
\stackrel{-\beta_l}{\longrightarrow} A_l$ 
which connects the fundamental alcove $A_\circ$ with its translation $A_{-\lambda}:=
A_\circ-\lambda.$ The operator $R^{[\lambda]}$ is defined as the composition 
\[ R^{[\lambda]}=(1+B_{\beta_l})\cdots (1+B_{\beta_1}), \] 
where $B_{\beta}$ is the {\em Bruhat operator} studied by Brenti, Fomin and Postnikov 
\cite{BFP}. 
Our construction of the model for $K_T(G/B)$ in $\B(V_W)$ 
is based on the operators $R^{[\lambda]}$. 
\medskip \\ 
{\bf Acknowledgements}\quad The second author would like to 
thank Anatol N. Kirillov for useful conversations. 

\section{Nichols-Woronowicz algebra}

The {\em Nichols-Woronowicz algebra} associated to a braided vector space 
is an analog of the polynomial ring in a braided setting. 
For details on the Nichols-Woronowicz algebra, see \cite{AS,Ba,Maj}. 

Let $(V,\psi)$ be a braided vector space, i.e. a vector space equipped with 
a linear isomorphism $\psi : V\otimes V \rightarrow V\otimes V$ such that 
the braid relation $\psi_i \psi_{i+1} \psi_i= \psi_{i+1} \psi_i \psi_{i+1}$ 
is satisfied on $V^{\otimes n}$; here 
$\psi_i$ is the linear endomorphism on $V^{\otimes n}$ obtained by applying 
$\psi$ on the $i$-th and $(i+1)$-st components. 
Throughout this paper, we assume that $V$ is a finite dimensional $\Q$-vector space. 
Let $w=s_{i_1}\ldots s_{i_l}$ be 
a reduced decomposition of an element $w\in S_n,$ where $s_i=(i,i+1)$ is
an adjacent transposition. Then the linear map $\Psi_w:= \psi_{i_1}\ldots \psi_{i_l}$ 
on $V^{\otimes n}$ is independent of the choice of a reduced decomposition of 
$w$ due to the braid relations. We define the {\em Woronowicz symmetrizer} on 
$V^{\otimes n}$ by $\sigma_n(\psi):= \sum_{w\in S_n} \Psi_w.$ 
Such a definition of the braided analog of the symmetrizer (or anti-symmetrizer) 
is due to Woronowicz \cite{Wo}. 

\begin{dfn} {\rm (see \cite{Ba} and \cite{Maj})} The Nichols-Woronowicz algebra 
$\B(V)$ associated to the braided vector space $(V,\psi)$ is defined as 
the quotient of the tensor algebra $T(V)$ by the ideal $\bigoplus_{n\geq 0}
{\rm Ker} (\sigma_n(\psi)).$ 
\end{dfn} 

\begin{rem} {\rm For a more systematic treatment, we need to work in a fixed 
braided category $\C$ of vector spaces. If the braided vector space 
$(V,\psi)$ is an object in the braided category $\C,$ the tensor algebra $T(V)$ 
has a natural braided Hopf algebra structure in $\C.$ It is known that 
the kernel $\bigoplus_{n\geq 0} {\rm Ker} (\sigma_n(\psi))$ is a Hopf ideal 
of $T(V).$ Hence, $\B(V)$ is also a braided Hopf algebra in $\C.$} 
\end{rem}

The following is the alternative definition of the Nichols-Woronowicz algebra
due to Andruskiewitsch and Schneider \cite{AS}. In \cite{AS}, the algebra 
$\B(V)$ is called the Nichols algebra. 

\begin{dfn} \cite{AS}
The graded braided Hopf algebra $\B(V)$ is called the Nichols-Woronowicz 
algebra if it satisfies the following conditions: 
\begin{enumerate}
\item[(1)] $\B_0(V)=\Q,$
\item[(2)] $V=\B_1(V)=\{ x\in \B(V) \; | \; \triangle (x)=x\otimes 1+ 1\otimes x \},$ 
\item[(3)] $\B(V)$ is generated by $\B_1(V)$ as a $\Q$-algebra. 
\end{enumerate}
\end{dfn} 

We use a particular braided vector space called the {\em Yetter-Drinfeld module} 
in the subsequent construction. 
Let $\Gamma$ be a finite group. 

\begin{dfn}
A $\Q$-vector space $V$ is called a Yetter-Drinfeld module over $\Gamma$ if 
\begin{enumerate} 
\item $V$ is a $\Gamma$-module, 
\item $V$ is $\Gamma$-graded, i.e. $V=\bigoplus_{g\in \Gamma}V_g,$ where $V_g$ is a 
linear subspace of $V,$ 
\item for $h\in \Gamma$ and $v\in V_g,$ we have $h(v)\in V_{hgh^{-1}}.$ 
\end{enumerate}
\end{dfn} 

The category $^{\Gamma}_{\Gamma}\YD$ of the Yetter-Drinfeld modules over 
a fixed group $\Gamma$ is naturally braided. The tensor product of the objects 
$V$ and $V'$ of $^{\Gamma}_{\Gamma}\YD$ 
is again a Yetter-Drinfeld module with the $\Gamma$-action 
$g(v\otimes w)=g(v)\otimes g(w)$ 
and the $\Gamma$-grading $(V\otimes V')_g= \bigoplus_{h,h'\in \Gamma, \; hh'=g}
V_h\otimes V'_{h'}.$ 
The braiding between $V$ and $V'$ is defined by $\psi_{V,V'}(v\otimes w)= 
g(w)\otimes v,$ for $v\in V_g$ and $w\in V'.$ 

Fix a Borel subgroup $B$ in a simple Lie group $G.$ 
Let $\Delta$ be the set of roots, and 
$\Delta_+$ the set of the positive roots corresponding to $B.$ 
We define a Yetter-Drinfeld module 
\[V_W:= \bigoplus_{\alpha \in \Delta}\Q \cdot [\alpha]/([\alpha]+[-\alpha]) \]
over the Weyl group $W.$ The $W$-action on $V_W$ is given by $w([\alpha])= [w(\alpha)].$ 
The $W$-degree of the symbol $[\alpha]$ is the reflection $s_{\alpha}.$ Note that 
$[\alpha]^2=0$ for all $\alpha \in \Delta,$ in the associated Nichols-Woronowicz 
algebra $\B(V_W).$ It is also easy to see that $[\alpha][\beta]=[\beta][\alpha]$ 
when $s_{\alpha}s_{\beta}=s_{\beta}s_{\alpha}.$ 
The following proposition can be shown by checking the quadratic relations in 
$\B(V_{S_n})$ via direct computation of the symmetrizer. 

\begin{prop} 
Fix the standard orthonormal basis $\{ e_1,\ldots , e_n \}$ of $\Q^n.$ 
Let $\Delta=\{ e_i-e_j \; | \; 1\leq i\ne j \leq n \}$ be the 
root system of type $A_{n-1}.$  Then there exists a surjective 
algebra homomorphism 
\[ \begin{array}{cccc} 
\eta : & \E_n & \rightarrow & \B(V_{S_n}) \\ [0.05in]
& [i,j] & \mapsto & [e_i-e_j] . 
\end{array} \]
\end{prop}

\begin{conj}
{\rm The algebra homomorphism $\eta$ is an isomorphism.}
\end{conj}

This conjecture is now confirmed up to $n=6.$ 

Take the standard realization of a root system of rank two with respect to 
an orthonormal basis $\{e_i\}$ as follows: \smallskip \\
\indent $(A_1 \times A_1)$ : $\Delta^{A_1 \times A_1}_+= \{ a_1=e_1, 
\; a_2=e_2 \},$ \smallskip \\ 
\indent $(A_2)$ : $\Delta^{A_2}_+= \{ a_1=e_1-e_2, \; a_2=e_1-e_3, 
\; a_3=e_2-e_3 \},$ \smallskip \\ 
\indent $(B_2)$ : $\Delta^{B_2}_+= \{ a_1=e_1-e_2, \; a_2=e_1, \; 
a_3=e_1+e_2, \; a_4=e_2 \},$ \smallskip \\ 
\indent $(C_2)$ : $\Delta^{C_2}_+= \{ a_1=e_1-e_2, \; a_2=2e_1, \; 
a_3=e_1+e_2, \; a_4=2e_2 \},$ \smallskip \\ 
\indent $(G_2)$ : $\Delta^{G_2}_+=
\{ a_1=e_1-e_2, \; a_2=e_1-2e_2+e_3, \; a_3=-e_2+e_3, \; 
a_4=-e_1-e_2+2e_3, \; a_5=-e_1+e_3, 
\; a_6=-2e_1+e_2+e_3 \}.$ \medskip \\ 
If the set $\Delta$ of roots contains a subset of the form 
$\Delta^X_+,$ where $X=A_1 \times A_1,$ $A_2,$ 
$B_2,$ $C_2,$ or $G_2,$ then one can check that the following relations are satisfied 
in the algebra $\B(V_W)$, respectively (see also \cite{KM3}). \smallskip \\ 
\indent $(A_1 \times A_1)$ : $[a_1][a_2]=[a_2][a_1],$ \smallskip \\ 
\indent $(A_2)$ : $[a_1][a_2]+[a_2][a_3]=[a_3][a_1],$ \smallskip \\ 
\indent $(B_2,C_2)$ : $[a_1][a_2]+[a_2][a_3]+[a_3][a_4]=[a_4][a_1],$ \smallskip \\ 
\hspace*{1.95cm} 
$[a_1][a_2][a_3][a_2]+[a_2][a_3][a_2][a_1]+[a_3][a_2][a_1][a_2]+[a_2][a_1][a_2][a_3]=0,$ 
\smallskip \\ 
\hspace*{1.95cm} 
$[a_2][a_3][a_4][a_3]+[a_3][a_4][a_3][a_2]+[a_4][a_3][a_2][a_3]+[a_3][a_2][a_3][a_4]=0,$ 
\smallskip \\ 
\hspace*{1.95cm} 
$[a_1][a_2][a_3][a_4]=[a_4][a_3][a_2][a_1],$ \smallskip \\ 
\indent $(G_2)$ : $[a_1][a_2]+[a_2][a_3]+[a_3][a_4]+[a_4][a_5]+[a_5][a_6]=[a_6][a_1],$ 
\begin{align*} 
& \!\!\!\!\!\!\!\!\!\!\!\!\!\!\!\!\!\!\!\!\!\!\!\!\!\!\!\!
\!\!\!\!\!\!\!\!\!\!\!\! [a_1][a_2][a_1][a_3]+[a_1][a_3][a_1][a_2]+
[a_1][a_3][a_2][a_3]= \\ & \!\!\!\!\!\!\!\!\!\!\!\!\!\!\!\!\!\!\!\!\!\!\!\!\!\!
\!\!\!\!\!\!\!\!\!\!\!\!\!\! =[a_2][a_1][a_3][a_1]+[a_3][a_2][a_3][a_1]+[a_3][a_1][a_2][a_1], \\
& \!\!\!\!\!\!\!\!\!\!\!\!\!\!\!\!\!\!\!\!\!\!\!\!\!\!\!\!\!\!\!\!\!\!\!\!
\!\!\!\! [a_6][a_5][a_6][a_4]+[a_6][a_4][a_6][a_5]+[a_6][a_4][a_5][a_4]= 
\\ & \!\!\!\!\!\!\!\!\!\!\!\!\!\!\!\!\!\!\!\!\!\!\!\!\!\!\!\!\!\!\!\!\!\!\!\!\!\!\!\! 
= [a_5][a_6][a_4][a_6]+[a_4][a_5][a_4][a_6]+[a_4][a_6][a_5][a_6], \\ 
& \!\!\!\!\!\!\!\!\!\!\!\!\!\!\!\!\!\!\!\!\!\!\!\!\!\!\!\!\!\!\!\!\!\!\!\!
\!\!\!\! [a_1][a_2][a_4][a_5]+[a_2][a_3][a_4][a_5]+[a_2][a_3][a_5][a_6]+
[a_5][a_3][a_4][a_2]= \\
& \!\!\!\!\!\!\!\!\!\!\!\!\!\!\!\!\!\!\!\!\!\!\!\!\!\!\!\!\!\!\!\!\!\!\!\!\!\!\!\! 
= [a_3][a_4][a_2][a_3]+[a_3][a_4][a_3][a_4]+
[a_4][a_5][a_3][a_4]+[a_6][a_4][a_3][a_1], \\ 
& \!\!\!\!\!\!\!\!\!\!\!\!\!\!\!\!\!\!\!\!\!\!\!\!\!\!\!\!\!\!\!\!\!\!\!\!
\!\!\!\! [a_5][a_4][a_2][a_1]+[a_5][a_4][a_3][a_2]+[a_6][a_5][a_3][a_2]+
[a_2][a_4][a_3][a_5]= \\
& \!\!\!\!\!\!\!\!\!\!\!\!\!\!\!\!\!\!\!\!\!\!\!\!\!\!\!\!\!\!\!\!\!\!\!\!
\!\!\!\! = [a_3][a_2][a_4][a_3]+[a_4][a_3][a_4][a_3]+
[a_4][a_3][a_5][a_4]+[a_1][a_3][a_4][a_6], \\ 
& \!\!\!\!\!\!\!\!\!\!\!\!\!\!\!\!\!\!\!\!\!\!\!\!\!\!\!\!\!\!\!\!\!\!\!\!
\!\!\!\! [a_1][a_2][a_3][a_4][a_5][a_6]=[a_6][a_5][a_4][a_3][a_2][a_1] . 
\end{align*}
The complete set of the independent defining relations for the algebra 
$\B(V_W)$ has not yet been determined in general. 
The relations listed above imply the following. 

\begin{prop}\label{prop12}
The elements $h_{\alpha}:=1+[\alpha],$ $\alpha\in \Delta,$ satisfy the 
Yang-Baxter equations, i.e., if $\Delta$ contains a subset $\Delta'$ of the form 
$\Delta^X_+,$ where $X=A_1 \times A_1,$ $A_2,$ $B_2,$ $C_2,$ or $G_2,$ 
then the elements $h_{\alpha},$ $\alpha \in \Delta',$ satisfy the following 
equations, respectively: \smallskip \\ 
\indent $(A_1 \times A_1)$ $:$ $h_{a_1}h_{a_2}=h_{a_2}h_{a_1},$ \smallskip \\ 
\indent $(A_2)$ $:$ $h_{a_1}h_{a_2}h_{a_3}=h_{a_3}h_{a_2}h_{a_1},$ \smallskip \\
\indent $(B_2, C_2)$ $:$ $h_{a_1}h_{a_2}h_{a_3}h_{a_4}=h_{a_4}h_{a_3}h_{a_2}h_{a_1},$ 
\smallskip \\
\indent $(G_2)$ $:$ $h_{a_1}h_{a_2}h_{a_3}h_{a_4}h_{a_5}h_{a_6}=
h_{a_6}h_{a_5}h_{a_4}h_{a_3}h_{a_2}h_{a_1}.$
\end{prop}

For a general braided vector space $V,$ the elements $v\in V$ act on the algebra 
$\B(V^*)$ as braided differential operators. In the subsequent construction, we 
use the braided differential operator $D_{\alpha}$ for a positive root 
$\alpha$, whose action on the algebra $\B(V_W)$ is determined by the following conditions: 
\begin{enumerate}
\item[(0)] $D_{\alpha}(t)=0,$ for $t\in \B_0(V_W)=\Q,$ 
\item[(1)] $D_{\alpha}([\beta])= \delta_{\alpha,\beta},$ for $\alpha,\beta \in 
\Delta_+,$ 
\item[(2)] $D_{\alpha}(FF')= D_{\alpha}(F)F'+s_{\alpha}(F)D_{\alpha}(F')$ for 
$F,F' \in \B(V_W).$ 
\end{enumerate} 
We set $D_{\alpha}:=-D_{-\alpha}$ if $\alpha$ is a negative root. 
The following is a key lemma in the proof of the main theorem. 

\begin{lem} \label{lem11} {\rm (see \cite[Proposition 2.4]{MS} and 
\cite[Criterion 3.2]{Ba})} We have
\[ \bigcap_{\alpha \in \Delta_+}{\rm Ker}(D_{\alpha})= \B_0(V_W)=\Q. \] 
\end{lem}

\section{Alcove path and multiplicative Chevalley elements}

In this section, we define a family of elements $\Xi^{\lambda},$ $\lambda \in P,$ 
in the Nichols-Woronowicz algebra $\B(V_W)$ based on the construction of 
the path operators due to the first author and Postnikov \cite{LP}. 

Let $W_{\af}$ be the affine Weyl group of the dual root system 
$\Delta^{\vee}:=\{ \alpha^{\vee}= 
2\alpha /\langle\alpha,\alpha\rangle \; | \; \alpha \in \Delta \}.$ 
The affine Weyl group $W_{\af}$ is generated by the affine reflections 
$s_{\alpha, k},$ $\alpha\in \Delta,$ $k\in \Z,$ with respect to the affine 
hyperplanes 
$H_{\alpha,k}:= \{ \lambda \in \h^* \; | \; \langle\lambda ,\alpha^{\vee}\rangle=k \}.$ 
The connected components of $\h^* \setminus \bigcup_{\alpha \in \Delta,k\in \Z} 
H_{\alpha,k}$ are called alcoves. The fundamental alcove $A_\circ$ is the alcove defined 
by the inequalities $0<\langle\lambda, \alpha^{\vee}\rangle < 1,$ for all 
$\alpha \in \Delta_+.$ We now recall some concepts from \cite{LP}. 

\begin{dfn} \cite{LP} $(1)$ A sequence $(A_0,\ldots,A_l)$ of alcoves $A_i$ is called 
an {\em alcove path} if $A_i$ and $A_{i+1}$ have a common wall and $A_i \not= A_{i+1}.$ 

$(2)$ An alcove path $(A_0,\ldots,A_l)$ is called reduced if the length $l$ 
of the path is minimal among all alcove paths connecting $A_0$ and $A_l.$ 

$(3)$ We use the symbol $A_i \stackrel{\beta}{\longrightarrow} A_{i+1}$ when 
$A_i$ and $A_{i+1}$ have a common wall of the form $H_{\beta,k}$ and the direction 
of the root $\beta$ is from $A_i$ to $A_{i+1}.$ 
\end{dfn} 

Let $\{ \alpha_1,\ldots, \alpha_r \} \subset \Delta_+$ be the set of the simple 
roots. 
Let $\omega_i$ be the fundamental weight corresponding to 
the simple root $\alpha_i,$ i.e. $\langle\omega_i,\alpha_j^{\vee}\rangle 
\delta_{i,j}.$ 
Take an alcove path $A_0 \stackrel{-\beta_1}{\longrightarrow} \cdots 
\stackrel{-\beta_l}{\longrightarrow} A_l$ connecting $A_0=A_\circ$ and 

$A_l=A_{-\lambda}:= A_\circ-\lambda.$ The sequence of roots $(\beta_1,\ldots ,\beta_l)$ 
appearing here is called a $\lambda$-{\em chain}. 

\begin{dfn}
We define the elements 
$\Xi^{[\lambda]}$ in $\B(V_W)$ for $\lambda \in P$ 
by the formula 
\[ \Xi^{[\lambda]}= h_{\beta_l} \ldots h_{\beta_1}. \]
We call the elements $\Xi_i:=\Xi^{\omega_i}$ the multiplicative 
Chevalley elements. 
\end{dfn} 

The element $\Xi^{[\lambda]}$ is independent of the choice of an 
alcove path from $A_\circ$ to $A_{-\lambda}.$ Indeed, the elements 
$h_{\alpha}$ satisfy the Yang-Baxter equations and $h_{\alpha}h_{-\alpha}=1,$ 
so the argument of 
\cite[Lemma 9.3]{LP}, which is implicitly given by Cherednik \cite{Ch}, 
is applicable to our case. 

We use the following results from \cite{LP}. 

\begin{lem} \label{lem21} {\rm (\cite[Lemmas 12.3 and 12.4]{LP})} 
Let $A_0 \stackrel{-\beta_1}{\longrightarrow} \cdots 
\stackrel{-\beta_l}{\longrightarrow} A_l$ be an alcove path from $A_0= A_\circ$ to 
$A_l=A_{-\lambda}.$ 

$(1)$ The sequence $(\alpha_i,s_{\alpha_i}(\beta_1),\ldots ,
s_{\alpha_i}(\beta_l), -\alpha_i)$ is an $s_{\alpha_i}(\lambda)$-chain for 
$i=1,\ldots,r.$ 

$(2)$ Asume that $\beta_j=\pm \alpha_i$ for some $1\leq j \leq l$ and 
$1\leq i \leq r.$ Denote by $s$ the reflection with respect to the common wall 
of $A_{l-j}$ and $A_{l-j+1}.$ Then the sequence 
$(\alpha_i,s_{\alpha_i}(\beta_1),\ldots , s_{\alpha_i}(\beta_{j-1}),\beta_{j+1}, 
\ldots , \beta_l)$ is an $s(\lambda)$-chain. 
\end{lem}

\begin{prop} {\rm (\cite[Proposition 12.2]{LP})}
We have $\Xi^{[\lambda]}\cdot \Xi^{[\lambda']}=\Xi^{[\lambda'+\lambda]},$ for 
$\lambda,\lambda' \in P.$ In particular, the elements $\Xi^{[\lambda]}$ and 
$\Xi^{[\lambda']}$ commute. 
\end{prop}

\section{Main result}

Consider the {\em Demazure operator} $\pi_i: \Z[P] \rightarrow \Z[P]$ 
corresponding to the simple root $\alpha_i$, which is defined by the formula 
\[ \pi_i(f) := 
\frac{ f-s_{\alpha_i}(f)}{e^{\alpha_i}-1} \,.\] 
The operator $\pi_i$ is characterized by the following conditions: 
\begin{enumerate}
\item $\pi_i(e^{\omega_j})=\delta_{i,j} e^{\omega_i-\alpha_i},$ 
\item $\pi_i(fg)=\pi_i(f)g+s_{\alpha_i}(f)\pi_i (g).$ 
\end{enumerate}
We have an algebra homomorphism 
\[ \begin{array}{cccc} 
\varphi : & \Z[P] & \rightarrow & \Z [\Xi_1,\ldots ,\Xi_r] \\ [0.05in]
& e^{\omega_i} & \mapsto & \Xi_i. 
\end{array} \] 

\begin{prop}\label{prop31}
Let $f$ be an element in $\Z[P].$ We have 
\[ \Pi_i(\varphi (f))= \varphi (\pi_if), \] 
where $\Pi_i=h_{\alpha_i}^{-1} \circ D_{[\alpha_i]}.$
\end{prop}

\begin{proof} It is enough to check that the operator $\Pi_i$ 
satisfies the following conditions: 
\begin{enumerate}
\item $\Pi_i(\Xi_j)=\delta_{i,j} \Xi^{[\omega_i-\alpha_i]},$ 
\item $\Pi_i(FF')=\Pi_i(F)F'+s_{\alpha_i}(F)\Pi_i (F'),$ for $F,F' \in 
\Z [\Xi_1,\ldots ,\Xi_r].$ 
\end{enumerate}
Let $t_i=t_{-\omega_j}\in W_{\af}$ be the translation by $-\omega_i.$ 
Since the hyperplane of the form $H_{\alpha_j,k},$ $j\not= i,$ does not separate
the alcoves $A_\circ$ and $t_i^{-1}(A_\circ),$ the roots $\pm \alpha_j,$ 
$j\not= i,$ 
cannot appear as a component of the $\omega_i$-chain 
$(\beta_1,\ldots, \beta_l)$ corresponding to 
a reduced alcove path $A_0 \stackrel{-\beta_1}{\longrightarrow} \cdots 
\stackrel{-\beta_l}{\longrightarrow} A_l$ from $A_0=A_\circ$ to 
$A_l=t_i(A_\circ)$ (see \cite[Chapter 4]{Hu}). 
Hence, we have $\Pi_i(\Xi_j)=0$ if $j\not= i.$ Based on 
Lemma \ref{lem21} (2), we also have $\Pi_i(\Xi_i)=\Xi^{[\omega_i-\alpha_i]},$ 
so condition (1) follows.

Pick an element $f\in \Z[P]$ such that $\varphi (f)=F.$ We have
\[ \Pi_i(FF')= h_{\alpha_i}^{-1}D_{[\alpha]}(F)F'+h_{\alpha_i}^{-1}\cdot
\varphi(s_{\alpha_i}(f)) \cdot 
h_{\alpha_i} \cdot h_{\alpha_i}^{-1}D_{[\alpha]}(F') \] 
by applying the twisted Leibniz rule for $D_{[\alpha]}.$ From 
Lemma \ref{lem21} (1), one can see that $h_{\alpha_i}^{-1}\cdot \varphi(s_{\alpha_i}(f)) 
\cdot h_{\alpha_i}= s_{\alpha_i}(F).$ So condition (2) is satisfied. 
\end{proof}

\begin{thm}
The subalgebra $\Z[\Xi_1,\ldots ,\Xi_r]$ 
generated by the multiplicative Chevalley elements 
in the Nichols-Woronowicz algebra $\B(V_W)$ 
is isomorphic to the $K$-ring $K(G/B).$ 
\end{thm} 

\begin{proof} Define the algebra homomorphism $\e: \Z[P] \rightarrow \Z$ 
by the assignment $e^{\lambda} \mapsto 1,$ for all $\lambda \in P.$ 
Let $I \subset \Z[P]$ be the ideal generated by the elements of 
the form $f-\e (f),$ $f\in \Z[P]^W.$ Then the $K$-ring $K(G/B)$ is 
isomorphic to the quotient algebra $\Z[P]/I.$ 
Pick an element $g\in I.$ We have $\Pi_i(\varphi(g))= \varphi(\pi_i(g))= 0$ for 
$i=1,\ldots ,r,$ by Proposition \ref{prop31}. Since $w \circ D_{[\alpha]} \circ w^{-1}
= D_{[w(\alpha)]},$ we obtain $D_{[\alpha]}(\varphi(g))=0,$ for all 
$\alpha \in \Delta_+,$ and therefore $\varphi(g)=0$ by Lemma \ref{lem11}.
If $g\not\in I,$ there exists an operator $\varpi$ on $\Z[P]$, written as 
a linear combination of the composites of the multiplication operators and 
the operators $\pi_i$, such that the constant term of $\varpi(g)$ is nonzero. 
We conclude that ${\rm Im}(\varphi) \cong \Z[P]/I \cong K(G/B).$ 
\end{proof}

\begin{rem} 
(1) The idea of the proof of the above theorem is used in \cite[Sections 5 and 6]{KM3} 
for the root systems of classical type and of type $G_2.$ 
The multiplicative Dunkl elements $\Theta_i:=\Xi^{[e_i]}$ corresponding to
the components of the orthonormal basis $\{e_i\}$ are used in \cite{KM3}. 
The multiplicative Dunkl elements in the Fomin-Kirillov quadratic algebra 
$\E_n$ appear in \cite{Le,LY,Ru}. 

(2) For an arbitrary parabolic subgroup $P \supset B,$ the $K$-ring $K(G/P)$ of 
the homogeneous space $G/P$ is a subalgebra of $K(G/B).$ Hence, the algebra 
$\B(V_W)$ also contains $K(G/P)$ as a commutative subalgebra. 
\end{rem}

Bazlov \cite{Ba} has proved that the subalgebra in $\B(V_W)$ generated by the 
elements $[\alpha]$ corresponding to the simple roots $\alpha$ is isomorphic
to the nil-Coxeter algebra 
\[ NC_W:= \Z \langle u_1,\ldots ,u_r \rangle / ( u_i^2, \; (u_iu_j)^{[m_{ij}/2]}u_i^{\nu_{ij}}-
(u_ju_i)^{[m_{ij}/2]}u_j^{\nu_{ij}}, \; i=1,\ldots, r), \] 
where $m_{ij}$ is the order of $s_{\alpha_i}s_{\alpha_j}$ in $W,$ 
$[m_{ij}/2]$ stands for the integer part of $m_{ij}/2,$ and 
$\nu_{ij}:= m_{ij}-2[m_{ij}/2].$ 
In our case, we can show the following. 

\begin{cor}
The subalgebra in ${\rm End}(\B(V_W))$ generated by the operators $\Pi_1,\ldots, 
\Pi_r$ is isomorphic to the nil-Hecke algebra 
\[ NH_W:= \Z \langle T_1,\ldots, T_r \rangle / ( T_i^2+T_i, \; 
(T_iT_j)^{[m_{ij}/2]}T_i^{\nu_{ij}}-
(T_jT_i)^{[m_{ij}/2]}T_j^{\nu_{ij}}, \; i=1,\ldots, r) \] 
via the map given by $T_i \mapsto \Pi_i.$ 
\end{cor}

\begin{proof} One can check that the operators $\Pi_i$ satisfy the defining 
relations of $NH_W$ by direct computations. Since the map given by $T_i \mapsto \pi_i$ 
defines a faithful representation of $NH_W$ on $\Z[P]/I,$ the subalgebra 
generated by $\Pi_i,$ $i=1,\ldots ,r$ is isomorphic to $NH_W.$ 
\end{proof}

\section{Model of the equivariant $K$-ring}
The results in the previous section are generalized to the case 
of the $T$-equivariant $K$-ring $K_T(G/B).$ 
Our construction of the model for $K_T(G/B)$ is also parallel to 
the approach in \cite{LP}.

Since the Nichols-Woronowicz algebra 
$\B(V_W)$ is a braided Hopf algebra in the category of the Yetter-Drinfeld 
modules over $W,$ it is $W$-graded. Denote by $w_x$ the $W$-degree of a $W$-homogeneous 
element $x\in \B(V_W).$ Let $h$ be the Coxeter number and $P':=h^{-1}\cdot P \subset \h^*.$ 
The Weyl group $W$ acts on the group algebra $\Z[P']=\bigoplus_{\lambda \in P'}X^{\lambda}$ 
by $w(X^{\lambda})=X^{w(\lambda)},$ $w\in W.$ 
The twist map 
\[ \begin{array}{cccc} 
c: & \B(V_W) \otimes \Z[P'] & \rightarrow & \Z[P'] \otimes \B(V_W) \\ [0.05in]
& x \otimes X & \mapsto & w_x(X) \otimes x 
\end{array} \] 
gives an associative multiplication map $m$ on $\B(V_W)\langle P' \rangle :=
\Z[P'] \otimes \B(V_W)$ as follows: 
\[ m: \Z[P'] \otimes \B(V_W) \otimes \Z[P'] \otimes \B(V_W) 
\arrl{1\otimes c \otimes 1} 
\Z[P'] \otimes \Z[P'] \otimes \B(V_W) \otimes \B(V_W) \] 
\[ \arrll{m_{\Z[P']} \otimes m_{\B}} \Z[P'] \otimes \B(V_W) , \] 
where $m_{\Z[P']}$ and $m_{\B}$ are the multiplication maps on 
the algebras $\Z[P']$ and $\B(V_W)$, respectively. 
The algebras $\Z[P']$ and $\B(V_W)$ are considered as subalgebras 
of $\B(V_W)\langle P' \rangle.$ We have the commutation relation 
\[ [\alpha] \cdot X^{\lambda} = X^{s_{\alpha}(\lambda)} \cdot 
[\alpha], \; \; \alpha \in \Delta, \; \lambda \in P', \] 
in $\B(V_W)\langle P' \rangle.$ 

The subalgebra $\Z[P]^W \subset \Z[P']$ is viewed as a subalgebra of 
the representation ring $R(T)$ of the maximal torus via the isomorphism 
\[ \begin{array}{cccc}
\iota : & \Z[P] & \rightarrow & R(T) \\ [0.05in]
& e^{\lambda} & \mapsto & \chi^{\lambda} . 
\end{array} \] 
Let us consider an $R(T)$-algebra 
$\B_T(V_W) := R(T) \otimes_{\Z[P]^W} \B(V_W)\langle P' \rangle .$ 
We introduce the elements 
\[ H_{\alpha}:=X^{\rho/h} \cdot (X^{\alpha/h}+[\alpha]) \cdot 
X^{-\rho/h}, \; \; \alpha \in \Delta, \; 
\rho := \frac{1}{2}\left(\sum_{\beta \in \Delta_+}\beta\right), \] 
in the algebra $\B_T(V_W).$ 
Since the argument in the proof of \cite[Theorem 10.1]{LP} 
is applicable to our case, Proposition \ref{prop12} implies the following. 

\begin{lem} \label{lem41}
{\rm The elements $H_{\alpha},$ $\alpha \in \Delta,$ satisfy the Yang-Baxter 
equations in the algebra $\B_T(V_W).$} 
\end{lem} 

Let $(\beta_1,\ldots,\beta_l)$ be a $\lambda$-chain 
for a weight $\lambda \in P.$ Define the element 
\[ \Xi^{[\lambda]}_{\rm eq}:= H_{\beta_l}\ldots H_{\beta_1} \] 
in $\B_T(V_W).$ The element $\Xi^{[\lambda]}_{\rm eq}$ is independent of 
the choice of the $\lambda$-chain from Lemma \ref{lem41}. We also have 
$\Xi^{[\lambda]}_{\rm eq} \cdot \Xi^{[\lambda']}_{\rm eq}= 
\Xi^{[\lambda+\lambda']}_{\rm eq}$ 
by \cite[Proposition 12.2]{LP}. 

The braided differential operators $D_{\alpha}$ are naturally extended 
as $R(T)$-linear operators on $\B_T(V_W)$, being determined by the following conditions: 
\begin{enumerate} 
\item[(0)] $D_{\alpha}(X)=0,$ for $X \in \Z[P'],$ 
\item[(1)] $D_{\alpha}([\beta])= \delta_{\alpha,\beta},$ for $\alpha,\beta \in 
\Delta_+,$ 
\item[(2)] $D_{\alpha}(FF')= D_{\alpha}(F)F'+s_{\alpha}(F)D_{\alpha}(F')$ for 
$F,F' \in \B(V_W)\langle P' \rangle.$ 
\end{enumerate}

\begin{lem} \label{lem42}
In the algebra $\B_T(V_W),$ we have 
\[ \bigcap_{\alpha \in \Delta_+}{\rm Ker}(D_{\alpha})=
R(T) \otimes_{\Z[P]^W} \Q[P']. \] 
\end{lem}

\begin{proof} 
This follows immediately from Lemma \ref{lem11}. 
\end{proof}

The operator $\pi_i$ defined in the previous section 
is extended $R(T)$-linearly to the group algebra $R(T)[P]$, being 
determined by the following conditions: 
\begin{enumerate} 
\item $\pi_i(e^{\omega_j})= \delta_{i,j} e^{\omega_i-\alpha_i},$ 
\item $\pi_i(fg)=\pi_i(f)g+s_{\alpha_i}(f)\pi_i (g).$ 
\end{enumerate} 
Here, the action of $W$ on $R(T)$ is assumed to be trivial. 

Note that $K_T(G/B)$ is isomorphic to the quotient algebra 
$R(T)[P]/J,$ where $J$ is the ideal generated by the elements 
of form $f-\iota(f),$ $f\in \Z[P]^W.$ 

\begin{thm}
The subalgebra $R(T)[\Xi^{[\lambda]}_{\rm eq},\; \lambda \in P]$ of $\B_T(V_W)$ 
is isomorphic to the $T$-equivariant $K$-ring $K_T(G/B).$ 
\end{thm} 

\begin{proof} 
Let us consider the homomorphism between $R(T)$-algebras 
\[ \begin{array}{cccc} 
\psi : & R(T)[P] & \rightarrow & R(T)[\Xi^{[\lambda]}_{\rm eq},\; \lambda 
\in P] \\ [0.05in]
& e^{\lambda} & \mapsto & \Xi^{[\lambda]}_{\rm eq} . 
\end{array} \] 
We can see that 
\[ X^{\rho/h}(X^{-\alpha_i/h}-[\alpha_i])X^{-s_{\alpha_i}(\rho)/h}
D_{\alpha_i}(\psi(f))= \psi(\pi_i(f)) , \; f\in R(T)[P], \] 
in the same manner as in the proof of Proposition \ref{prop31}. 
Therefore, if $f\in \Z[P]^W,$ then $D_{\alpha}(\psi(f))=0,$ for all $\alpha 
\in \Delta_+.$ Based on Lemma \ref{lem42}, we have 
$\psi(f) \in R(T) \otimes_{\Z[P]^W} \Z[P'].$ Here, 
the constant term of $\psi(f)$ for $f\in \Z[P]^W$ is in $\Z[P]^W$, and 
therefore equals $\iota(f)$ (see also \cite[Proposition 14.5]{LP}). So 
we obtain $\psi(f)=0$ for $f\in J.$ 
On the other hand, if $f \not\in J,$ there exists an operator 
$\varpi$ on $R(T)[P]$, written as a linear combination 
of the composites of the multiplication operators and the operators $\pi_i$, 
such that the constant term of $\varpi(f)$ is nonzero. 
We conclude that ${\rm Ker}(\psi)=J.$ 
\end{proof}

\end{document}